\documentclass{amsart}

\usepackage{amsmath,amssymb,amsfonts, amscd,txfonts}
\newcommand{\Hom}{\normalfont\mbox{Hom}\,}
\newcommand{\Ext}{\normalfont\mbox{Ext}\,}
\newcommand{\Tor}{\normalfont\mbox{Tor}\,}
\newcommand{\pd}{\normalfont\mbox{pd}\,}
\newcommand{\id}{\normalfont\mbox{id}\,}
\newcommand{\fd}{\normalfont\mbox{fd}\,}
\newcommand{\gd}{\mbox{Gdim}\,}
\newcommand{\wdim}{\normalfont\mbox{wdim}\,}

\newcommand{\gldim}{\normalfont\mbox{gldim}\,}
\newcommand{\wGgldim}{\normalfont\mbox{wGgldim}\,}
\newcommand{\Ggldim}{\normalfont\mbox{Ggldim}\,}
\newcommand{\Gpd}{\normalfont\mbox{Gpd}\,}
\newcommand{\Gid}{\normalfont\mbox{Gid}\,}
\newcommand{\Gfd}{\normalfont\mbox{Gfd}\,}
\newcommand{\im}{\normalfont\mbox{Im}\,}

\newcommand{\FPD}{\normalfont\mbox{FPD}\,}
\newcommand{\cd}{\normalfont\mbox{cd}\,}
\newcommand{\cotd}{\normalfont\mbox{cot.D}\,}

\theoremstyle{plain}
\newtheorem{theorem}{Theorem}[section]
\newtheorem{lemma}[theorem]{Lemma}

\theoremstyle{definition}
\newtheorem{definition}[theorem]{Definition}

\newtheorem{example}[theorem]{Example}

\theoremstyle{Definition and Notation}

\catcode`\ç=13
\defç{\c{c}}
\catcode`\é=13
\defé{\'e}
\catcode`\à=13
\defà{\`a}
\catcode`\è=13
\defè{\`e}
\catcode`\â=13
\defâ{\^a}
\catcode`\ù=13
\defù{\`u}
\catcode`\ê=13
\defê{\^e}
\catcode`\î=13
\defî{\^\i}
\catcode`\ô=13
\defô{\^o}
\begin{document}

\title[Gorenstein global dimension of an amalgamated...]{Gorenstein global dimension of an amalgamated duplication of a coherent ring along an ideal}

\author[N. Mahdou]{Najib Mahdou}
\address{Department of Mathematics, Faculty of Science and Technology of Fez, Box 2202, University S. M.
Ben Abdellah Fez, Morocco}
 \email{mahdou@hotmail.com}

\author[M. Tamekkante]{Mohammed Tamekkante}
\address{Department of Mathematics, Faculty of Science and Technology of Fez, Box 2202, University S. M.
Ben Abdellah Fez, Morocco}
 \email{tamekkante@yahoo.fr}





\begin{abstract} In this paper, we study the Gorenstein global dimension
of an \emph{amalgamated duplication} of a coherent ring along a
regular principal ideal.
\end{abstract}

\maketitle

\section{Introduction} Throughout this paper, all rings are commutative with identity
element, and all modules are unital. If $M$ is an $R$-module, we use
$\pd_R(M)$, $\id_R(M)$ and $\fd_R(M)$ to denote, respectively, the
classical projective, injective and flat dimensions of $M$.  We use also $\gldim(R)$ and $\wdim(R)$ to denote, respectively, the classical global and weak dimension of $R$.\\
We say that an ideal is principal regular if it is generated by a
regular element, i.e., a
non-zerodivisor element. \\
\emph{The amalgamated duplication of a ring $R$ along  an
$R$-submodule of the total ring of quotients $T(R)$}, introduced by
D'Anna and Fontana and denoted by $R\bowtie E$ (see \cite{D'Anna
0,D'Anna,D'Anna basic}), is  the following subring of $R\times T(R)$
(endowed with the usual componentwise operations):
$$R\bowtie E:=\{(r,r+e)\mid  r\in R\ \text{ and } e\in E\}.$$

It is obvious that, if in the $R$-module $R\oplus E$ we introduce a
multiplicative structure by setting $(r,e)(s, f):=(rs, rf+se+ef)$,
where $r,s\in R$ and $e,f\in E$, then we get the ring isomorphism
$R\bowtie E\cong R\oplus E$. When $E^2=0$, this new construction
coincides with the \emph{Nagata's idealization}. One main difference
between  this  constructions, with respect to the idealization (or
with respect to any commutative extension, in the sense of Fossum)
is that the ring $R\bowtie E$ can be a reduced ring and it is always
reduced if R is a domain (see \cite{D'Anna 0,D'Anna basic}). If
$E=I$ is an ideal in $R$, then the ring $R\bowtie I$ is a subring of
$R\times R$. This extension has been studied, in the general case,
and from the different point of view of pullbacks, by D'Anna and
Fontana \cite{D'Anna basic}. As it happens for the idealization, one
interesting application of this construction is the fact that it
allows to produce rings satisfying (or not satisfying) preassigned
conditions. Recently, D'Anna proved that, if $R$ is a local
Cohen-Macaulay ring with canonical module $\omega_R$, then $R\bowtie
I$ is a Gorenstein ring if and only if $I\cong \omega_R$ (see
\cite{D'Anna 0}). Note also that this construction has already been
applied, by Maimani and Yassemi for studying questions concerning
the diameter and girth of the zero-divisor graph of a ring (see
\cite{Yassemi}).\\
Recently in \cite{chhitimahdou}, the authors study some homological
properties and
coherence of the \emph{amalgamated duplication} of a ring along an ideal.\\
In this paper we study some particular cases of the
\emph{amalgamated duplication} extension from  the Gorenstein
homological point of view.

For a two-sided Noetherian ring $R$, Auslander and Bridger
\cite{Aus bri} introduced the $G$-dimension, $\gd_R (M)$, for
every finitely generated $R$-module $M$. They showed that there is
an inequality $\gd_R (M)\leq \pd_R (M)$ for all finite $R$-modules
$M$, and equality holds if $\pd_R (M)$ is finite.

Several decades later, Enochs and Jenda \cite{Enochs,Enochs2}
defined the notion of Gorenstein projective dimension
($G$-projective dimension for short), as an extension of
$G$-dimension to modules which are not necessarily finitely
generated, and the Gorenstein injective dimension ($G$-injective
dimension for short) as a dual notion of Gorenstein projective
dimension. Then, to complete the analogy with the classical
homological dimension, Enochs, Jenda and Torrecillas \cite{Eno
Jenda Torrecillas} introduced the Gorenstein flat dimension. Some
references are
 \cite{Christensen, Christensen
and Frankild, Enochs, Enochs2, Eno Jenda Torrecillas, Holm}.\\

Recall that an $R$-module $M$ is called
Gorenstein projective if, there exists an exact sequence of
projective $R$-modules:
$$\mathbf{P}:...\rightarrow P_1\rightarrow P_0\rightarrow
P^0\rightarrow P^1\rightarrow ...$$ such that $M\cong
\im(P_0\rightarrow P^0)$ and such that the operator $\Hom_R(-,Q)$
leaves $\mathbf{P}$ exact whenever $Q$ is
projective. The resolution $\mathbf{P}$ is called a complete projective resolution. \\
The Gorenstein injective $R$-module are defined
dually.\\
 And an $R$-module $M$ is called l
Gorenstein flat if, there exists an exact sequence of flat $R$-modules:
$$\mathbf{F}:...\rightarrow F_1\rightarrow F_0\rightarrow
F^0\rightarrow F^1\rightarrow ...$$ such that $M\cong
\im(P_0\rightarrow P^0)$ and such that the operator $I\otimes_R-$
leaves $F$ exact whenever $I$ is a injective. The resolution $\mathbf{F}$ is called complete flat resolution.\\

The  Gorenstein projective, injective and flat
dimensions are defined in term of resolution and denoted by $\Gpd(-)$, $\Gid(-)$ and $\Gfd(-)$ respectively (see \cite{Christensen, Enocks and janda, Holm}).\\

In \cite{Bennis and Mahdou2}, the authors prove the equality:
$$\sup\{\Gpd_R(M)\mid \text{M is an R-module}\}=\sup\{\Gid_R(M)\mid \text{M is an R-module}\}$$
They called the common value of the above quantities the
Gorenstein global dimension of $R$ and denoted it by
$\Ggldim(R)$. Similarly, they set
$$\wGgldim(R)=\sup\{\Gfd_R(M)\mid \text{M is an R-module}\}$$
which they called the  weak Gorenstein global dimension of
$R$. \\
In this paper we study the Gorenstein homological
dimension of the class of rings with the form $R\bowtie xR$ where
$R$ is coherent and $x$ is a  nonunit regular element of $R$.\\
 As shown in \cite[Corollary 1.2]{Bennis and Mahdou2}, the Gorenstein homological dimensions of rings are
refinement of the classical homological ones with equality if the
weak dimension is finite. Hence, the importance of these dimensions
appears when the weak dimension is infinite. Thus, before
continuing, it is naturel (as motivation)  to show what happens
with the $\wdim(R\bowtie xR)$.\\
Consider the  short exact sequences of $R\bowtie xR$-modules:
$$(\star)\quad 0\longrightarrow 0\times xR\stackrel{\iota}\longrightarrow R\bowtie
xR\stackrel{\varepsilon}\longrightarrow R \longrightarrow 0$$ where
$i$ is the injection and $\varepsilon(r,r+r'x)=r$, and
$$(\star\star) \quad 0 \longrightarrow R\stackrel{\mu}\longrightarrow
R\bowtie xR \stackrel{\nu}\longrightarrow 0\times xR\longrightarrow
0$$ where $\mu(r)=(rx,0)$ and $\nu(r,r+xr')=(0,(r+r'x)x)$\\
It is easy, from $(\star)$ and $(\star\star)$, to see that  $R$ and
$0\times xR$ are finitely presented $R\bowtie xR$-modules since
$0\times R=(0,x)R\bowtie xR$ and $R\cong R\bowtie xR/(0\times xR)$
are finitely generated. And, for every $R\bowtie xR$-module $M$,
\begin{itemize}
  \item $\Tor_{R\bowtie xR}(M,R)=\Tor_{R\bowtie xR}^2(M,0\times xR)=\Tor_{R\bowtie
  xR}^3(M,R)=...$.
  \item $\Tor_{R\bowtie xR}(M,0\times xR)=\Tor^2_{R\bowtie xR}(M,R)=\Tor^n_{R\bowtie xR}(M,0\times
  xR)=...$.
\end{itemize}

Thus, if $\wdim(R\bowtie xR)$ is finite, then $R$ and $0\times xR$ are flat
(then projective since they are finitely presented). But $0\times
xR=(0,x)R\bowtie xR$ can not be a direct summand of $R\bowtie xR$.
Deny. By \cite[Theorem 1.4.5]{Glaz}, there is an element
$(r,r+xr')\in R\bowtie xR$ such that $(0,x^2)(r,r+xr')=(0,x)$. That
means that $x(r+xr')=1$. Then $x$ is invertible, a contradiction. Then, we
conclude that $wdim(R\bowtie xR)=\infty$. So, it is interesting to
study the Gorenstein global and weak dimensions of
$R\bowtie xR$.\\

The main result of this paper is stated as follows:

\begin{theorem}\label{main result2}

Let $R$ be a coherent ring which contains a nonunit regular  element
$x$. Then,
\begin{enumerate}
  \item $\wGgldim(R\bowtie xR)=\wGgldim(R)$, and
  \item $\Ggldim(R\bowtie xR)=\Ggldim(R)$.
\end{enumerate}

\end{theorem}

\section{Main results}

To prove this Theorem  we need to introduce some other notions (the
 $n$-FC rings and the  global cotorsion dimension of
rings) and we need also a several Lemmas.\\

\begin{lemma}\label{id} Let $R$ be a ring and $I$ be a proper ideal of $R$ such that the canonical
(multiplication) map $R \rightarrow \Hom_R(I, I)$ is an isomorphism
and $\Ext^i_R(I,I)=0$ for all $i>0$. Then, $\id_{R\bowtie I}(R\bowtie
I)=\id_R(I)$.\\
In particular, if $x$ is a regular element of $R$ then,
$\id_{R\bowtie xR}(R\bowtie xR)=\id_R(R)$.
\end{lemma}

\begin{proof} Consider a $k$-step injective resolution of $I$ as an
$R$-module:
$$(\star)\quad 0\rightarrow I\rightarrow I_1\rightarrow
...\rightarrow I_k \rightarrow K \rightarrow 0.$$ By \cite[Theorem
2.1]{D'Anna}, $R\bowtie I\cong_{R\bowtie I} \Hom_R(R\bowtie I,I)$.
On the other hand, by hypothesis, $\Ext_R^i(R\bowtie I, I)\cong
\Ext_R^i(R\oplus I,I)=0$ for all $i>0$. Then, if we apply
$\Hom_R(R\bowtie xR,-)$ to $(\star)$, we deduce a $k$-step
injective resolution of $R\bowtie I$ as an $R\bowtie I$-module
(see that $\Hom_R(R\bowtie I,I_i)$ is an injective $R\bowtie
I$-module for each $1\leq i\leq k$ since $I_i$ is an injective
$R$-module):
$$0\rightarrow R\bowtie I\rightarrow \Hom_R(R\bowtie
I,I_1)\rightarrow ... \rightarrow \Hom_R(R\bowtie I,I_k)\rightarrow
\Hom_R(R\bowtie I,K)\rightarrow 0$$

If $K$ is an injective $R$-module, then $\Hom_R(R\bowtie I,K)$ is an
injective $R\bowtie I$-module. Hence, $\id_{R\bowtie I}(R\bowtie
I)\leq \id_R(I)$.\\ On the other hand,
\begin{center}
$\Hom_{R\bowtie I}(R,\Hom_R(R\bowtie I, K))\cong_R \Hom_R(R,K)\cong_RK$
(\cite[Lemma 11.61(iv)]{Rotman})
\end{center}
 (we use the ring morphism
$R\rightarrow R\bowtie I$ defined by $r\mapsto (r,r)$). Then, if
$\Hom_R(R\bowtie I,K)$ is an injective $R\bowtie I$-module, $K$ is
also an injective $R$-module. Hence, $\id_R(I)\leq \id_{R\bowtie
I}(R\bowtie I)$.\\ Consequently, $\id_R(I)=\id_{R\bowtie I}(R\bowtie
I)$, as desired.\\
Now suppose that $x$ is a regular element of $R$ and consider the
multiplication morphism of $R$-modules $\varepsilon:R\rightarrow
\Hom_R(xR,xR)$ defined by, for every $r,a\in R$,
$\varepsilon(r)(xa)=rxa$. We claim that $\varepsilon$ is an
isomorphism. If $\varepsilon(r)=0$ we have $\varepsilon(r)(x)=rx=0$.
Then, $r=0$ since $x$ is a regular element of $R$. Now, let
$\varphi\in \Hom_R(xR,xR)$. Since $x$ is regular there is a
 unique $a\in R$ such that $\varphi(x)=ax$. Thus,  for every $y\in R$ we have
$\varphi(xy)=y\varphi(x)=yax=a(xy)$. Then, $\varphi=\varepsilon(a)$
and then $\varepsilon$ is an isomorphism. Consequently, the desired
particular result follows immediately from the first part of this
proof.
\end{proof}

\begin{definition}(\cite{Stenstrom} and \cite{Garkusha}) Let $R$ be a
ring and $M$ be an  $R$-module.
\begin{enumerate}
    \item The $FP$-injective dimension of $M$, denoted by $FP$-$\id_R(M)$, is the least positive integer $n$ such that $\Ext^{n+1}_R(P,M) = 0$ for every finitely
     presented $R$-module $P$.
    \item A ring $R$ is said to be $n$-$FC$, if it is coherent and it has
    self-$FP$-injective at most at $n$ (i.e., $FP$-$\id_R(R)\leq n$).
\end{enumerate}
\end{definition}
In the Definition above, if $n=0$ the ring is called an $FC$ ring.
Using \cite[Theorems 6 and 7]{Chen}, we deduce the following Lemma.
\begin{lemma}\label{n-FC}
Let $R$ be a coherent ring. The following statements are equivalent for an
integer $n\geq 0$.
\begin{enumerate}
    \item $R$ is $n-FC$.
    \item $\wGgldim(R)\leq n$
\end{enumerate}
\end{lemma}

\begin{lemma}\label{FP-id}
Let $R$ be a ring which contains a regular element $x$. Then,
$FP$-$\id_R(R)=FP$-$\id_{R\bowtie xR}(R\bowtie xR)$.
\end{lemma}
\begin{proof}
First, we claim that $FP$-$\id_R(R)\leq FP$-$\id_{R\bowtie
xR}(R\bowtie xR)$. Obviously, if $FP$-$\id_{R\bowtie xR}(R\bowtie
xR)=\infty$, then the inequality holds. Otherwise, we may assume
that $FP$-$\id_{R\bowtie xR}(R\bowtie xR)= n<\infty$. Let $M$ be a
finitely presented $R\bowtie xR$-module. From the short exact
sequence of $R\bowtie xR$-modules:
$$(\star)\quad 0\rightarrow 0\times
xR\stackrel{\iota}\rightarrow R\bowtie
xR\stackrel{\varepsilon}\rightarrow R \rightarrow 0$$ where
$\iota$ is the injection and $\varepsilon(r,r+r'x)=r$, we see
that $M$ is also a finitely presented $R\bowtie xR$-module (via
$\varepsilon$ and by using \cite[Theorem 2.1.8(2)]{Glaz} since
$0\times xR=(0,x)R\bowtie xR$ is finitely generated). Applying
\cite[Theorem 11.66]{Rotman} to the diagonal embedding $\varphi:R
\rightarrow R\bowtie xR$ defined by $\varphi(r)=(r,r)$ (see that
$R\bowtie xR$ is a free $R$-module), we obtain:
$$(\ast)\qquad \Ext_R^{n+1}(M,xR)\cong \Ext_{R\bowtie xR}^{n+1}(M,\Hom_R(R\bowtie
xR,xR)).$$ But $xR\cong_RR$ and  $R\bowtie xR\cong_{R\bowtie
xR}\Hom_R(R\bowtie xR,xR)$ (see the proof of Lemma \ref{id}).
Then, $(\ast)$ becomes:  $$\Ext_R^{n+1}(M,R)\cong \Ext_{R\bowtie
xR}^{n+1}(M,R\bowtie xR)=0.$$ Thus,
$FP$-$\id_R(R)\leq n$, as desired.\\
Secondly, we claim that $FP$-$\id_{R\bowtie xR}(R\bowtie xR)\leq
FP$-$\id_R(R)$. The case $FP$-$\id_R(R)=\infty$ is obvious. So, we
may assume  that $FP$-$\id(R)= n<\infty$. Let $M$ be  a finitely
presented $R\bowtie xR$-module. Since $R\bowtie xR$ is a finitely
generated free $R$-module via the morphism of rings $\varphi$
($R\bowtie xR\cong_RR^2$), $M$ is also a finitely presented
$R$-module. Hence,
 by \cite[Theorem 11.66]{Rotman}, we have:
$$\Ext_{R\bowtie xR}^{n+1}(M,R\bowtie xR)\cong \Ext_{R\bowtie
xR}^{n+1}(M,\Hom_R(R\bowtie xR,xR)) \cong \Ext_R^{n+1}(M,xR)=0$$ since
$FP$-$\id_R(xR)=FP$-$\id_R(R)= n$. It follows that $FP$-$\id_{R\bowtie
xR}(R\bowtie xR)\leq n$. Consequently, $FP$-$\id_R(R)=FP$-$\id_{R\bowtie
xR}(R\bowtie xR)$, as desired.
\end{proof}

\begin{lemma}\label{passage}
Let $R$ be a ring which contains a regular element $x$ and $M$ be an
$R\bowtie xR$-module. If $\Tor_{R\bowtie xR}^i(M,R)=0$ for every $i>0$,
then $\pd_{R\bowtie xR}(M)=\pd_R(M\otimes_{R\bowtie xR}R)$.\\
 In particular, $M$ is a
projective $R\bowtie xR$-module if and only if, $M\otimes_{R\bowtie
xR}R$ is a projective $R$-module.
\end{lemma}

\begin{proof} Recall that every $R\bowtie xR$-module $M$ (resp.,
every $R\bowtie xR$-morphism $f$) is an $R$-module (resp. an
$R$-morphism) via the diagonal embedding of  rings
$\upsilon:R\rightarrow R\bowtie xR$ defined by $r\mapsto (r,r)$.
Explicitly, for every $r,r' \in R$ and $m\in M$, we have the
$R$-modulation
$r.m:=(r,r)m$.\\
Clearly, we have $\pd_R(M\otimes_{R\bowtie xR}R)\leq \pd_{R\bowtie
xR}(M)$ since $\Tor_{R\bowtie xR}^i(M,R)=0$. Hence, we have to prove
the other inequality.\\
 The proof will be by induction on
$n:=\pd_R(M\otimes_{R\bowtie xR}R)$.\\
 First, suppose that
$M\otimes_{R\bowtie xR}R$ is projective and consider the short exact
sequence of $R\bowtie xR$-modules:
$$(\star)\quad 0\longrightarrow 0\times xR\stackrel{\iota}\longrightarrow R\bowtie
xR\stackrel{\varepsilon}\longrightarrow R \longrightarrow 0$$ where
$i$ is the injection and $\varepsilon(r,r+r'x)=r$ for every $r,r'\in R$.\\
Since $\Tor_{R\bowtie xR}(M,R)=0$, if we tensor $(\star)$ with
$M\otimes_{R\bowtie xR}-$ we obtain the short exact sequence of
$R\bowtie xR$-modules:
$$(\star\star)\quad 0\longrightarrow M\otimes_{R\bowtie xR}(0\times xR)\stackrel{1_M\otimes\iota}\longrightarrow M\otimes_{R\bowtie xR}R\bowtie
xR\stackrel{1_M\otimes\varepsilon}\longrightarrow M\otimes_{R\bowtie
xR}R \longrightarrow 0$$ It is also an exact sequence of $R$-modules
via $\upsilon$.\\
On the other hand, consider the $R$-module $0\times xR$ (the
modulation of $R$ over $0\times xR$ is the naturel one defined by
$r(0,xa)=(0,rxa)$). Since $xR\cong_R0\times xR$ naturally, we have
the naturel isomorphism of $R$-modules $$(1)\qquad
(M\otimes_{R\bowtie xR}R)\otimes_R(0\times
xR)\cong_R(M\otimes_{R\bowtie xR}R)\otimes_RxR$$ Moreover, if we
consider the ring map $\varphi:R\bowtie xR\rightarrow R$ defined by
$(r,r+xr')\mapsto r+xr'$, using \cite[Lemma 11.61(i)]{Rotman}, we
have the naturel isomorphism of $R$-modules:
$$(2)\qquad (M\otimes_{R\bowtie xR}R)\otimes_R(0\times
xR)\cong_R M\otimes_{R\bowtie xR}(0\times xR)$$ and it is clear that
the  modulation of $R\bowtie xR$ over $0\times xR$  via $\varphi$ is
defined by setting:
$$(r,r+xr')\cdot(0,xa)=\varphi((r,r+xr')(0,xa)=(r+xr')(0,xa)=(0,(r+xr')xa)$$
Clearly, this modulation is the same as the naturel modulation of
$R\bowtie xR$ over the ideal $0\times xR$. If we denote
$\overline{M}:=M/(0\times xR)M$, we have the naturel isomorphisms
of $R$-modules:  $M\otimes_{R\bowtie xR}R\cong_{\phi_1}
\overline{M}$ and $M\otimes_{R\bowtie xR}R\bowtie xR\cong_{\phi_2}
M$ . Thus, from $(1)$, $(2)$ and $\phi_1$, we have the naturel
isomorphism of $R$-modules:

$$\begin{array}{cccc}
  \phi_3: & M\otimes_{R\bowtie xR}(0\times xR) & \stackrel{\simeq}\rightarrow & \overline{M}\otimes_RxR \\
 & m\otimes (0,xr) & \mapsto & \overline{m}\otimes xr
\end{array}$$

 It is not hard to
check that  the following  diagram of $R$-modules is commutative:

$$\begin{array}{ccccccc}
  0\longrightarrow  & M\otimes_{R\bowtie xR}(0\times xR) & \stackrel{1_M\otimes\iota}\longrightarrow  & M\otimes_{R\bowtie xR}R\bowtie
xR & \stackrel{1_M\otimes\varepsilon}\longrightarrow&
M\otimes_{R\bowtie xR}
R  &\longrightarrow 0 \\
   & \phi_3\parallel \wr&  & \phi_2\parallel \wr&  & \phi_1\parallel \wr& \\
   0\longrightarrow & \overline{M}\otimes_RxR & \stackrel{\alpha}\longrightarrow  & M & \stackrel{\pi}\longrightarrow  & \overline{M} &
   \longrightarrow 0
\end{array}$$

where $\alpha$ and $\pi$ are defined as follows:
$\alpha(\overline{m}\otimes xr)=(0,xr)m$,  for every $m\in M$ and
$r\in R$   and $\pi(m)=\overline{m}$ is the canonical projection.
Hence, since the upper sequence is exact, so is the down one.
Moreover, this sequence splits since $\overline{M}$ is a
projective $R$-module and there is an $R$-morphism
$\pi':\overline{M}\rightarrow M$ such that $\pi\circ
\pi'=id(\overline{M})$. Hence, we deduce the following isomorphism
of $R$-modules:
$$\begin{array}{cccc}
  \beta: & \overline{M}\oplus (\overline{M}\otimes_RxR)& \rightarrow &M \\
   & (\overline{m},\overline{m'}\otimes x) & \mapsto & \pi'(\overline{m})+\alpha(\overline{m'}\otimes
   x)=\pi'(\overline{m})+(0,x)m'
\end{array}$$

Give the left hand side of this isomorphism an $R\bowtie xR$-modules
structure by setting:
$$(r,r+xr')\ast(\overline{m},\overline{m'}\otimes x)=(r\overline{m},\overline{m}\otimes
r'x+\overline{m'}\otimes (r+r'x)x)$$ (it is easy to check that the
structure above is a modulation and we denote the $R\bowtie
xR$-module obtaining by $\widehat{M}$).\\
 Seen that  $\overline{\pi'(\overline{m})}=\pi\circ \pi'(\overline{m})=\overline{m}$ we have:

\begin{eqnarray*}
  \beta((r,r+xr')\ast(\overline{m},\overline{m'}\otimes x))&=&\beta(r\overline{m},\overline{m}\otimes
r'x+\overline{m'}\otimes(r+r'x)x )\\ &=&
\beta(r\overline{m},\overline{\pi'(\overline{m})}\otimes
r'x+\overline{m'}\otimes(r+r'x)x
)\\ &=& r.\pi'(\overline{m})+(0,r'x)\pi'(\overline{m})+(0,(r+r'x)x)m' \\
&=& (r,r)\pi'(\overline{m})+(0,r'x)\pi'(\overline{m})+(0,(r+r'x)x)m' \\
  &=& (r,r+xr')(\pi'(\overline{m})+(0,x)m') \\
 &=&(r,r+r'x)\beta(\overline{m},\overline{m'}\otimes
xr)
\end{eqnarray*}
 Thus, $\beta$ is an isomorphism of $R\bowtie xR$-modules from $\widehat{M}$ into $M$.\\
Recall that $\overline{M}$ is a projective $R$-module and consider
an $R$-morphisms: $\rho(:=(\rho_i)_{i\in I}):
\overline{M}\rightarrow R^{(I)}$ and $\sigma: R^{(I)}\rightarrow
\overline{M}$ for a free $R$-module $R^{(I)}$ such that
$\varrho\circ \rho =id(\overline{M})$.\\
Using $\rho$ and $\sigma$, we construct  the $R$-morphisms:
$$\begin{array}{cccc}
  \widetilde{\rho}: & \overline{M}\oplus (\overline{M}\otimes_RxR)& \rightarrow & (R\bowtie xR)^{(I)} \\
   & (\overline{m},\overline{m'}\otimes x) & \mapsto &
(\rho_i(\overline{m}),\rho_i(\overline{m})+x\rho_i(\overline{m'}))_{i\in
I}
\end{array}$$
and
$$\begin{array}{cccc}
  \widehat{\varrho}: & (R\bowtie xR)^{I} & \rightarrow  &\overline{M}\oplus (\overline{M}\otimes_RxR) \\
   & ((r_i,r_i+r'_ix))_{i\in I} & \rightarrow & (\varrho((r_i)_{i\in I}),\varrho((r'_i)_{i\in I})\otimes x)
\end{array}$$

It is clear that $\widetilde{\rho}$ and $\widehat{\varrho}$ are well
defined since $x$ is regular. Moreover, we have:
\begin{eqnarray*}
\widetilde{\rho}((r,r+xr')\ast(\overline{m},\overline{m'}\otimes
x))&=&\widetilde{\rho}(r\overline{m},\overline{m}\otimes
r'x+\overline{m'}\otimes(r+r'x)x )\\ & = & (r\rho_i(\overline{m}),r\rho_i(\overline{m})+r'x\rho_i(\overline{m})+(r+xr')x\rho_i(\overline{m'}))_{i\in I} \\
  & = & (r,r+xr')(\rho_i(\overline{m}),\rho_i(\overline{m})+x\rho_i(\overline{m}))_{i\in
I} \\
& = & (r,r+xr')\widetilde{\rho}(\overline{m},\overline{m'}\otimes x)
\end{eqnarray*}
and
\begin{eqnarray*}
  \widetilde{\varrho}((r,r+xr')((r_i,r_i+r'_ix))_{i\in I}) & = &  \widetilde{\varrho}(((rr_i,rr_i+(rr'_i+r'r_i+xr'r'_i)x))_{i\in I}) \\
   & = &(r\varrho((r_i)_{i\in I}),r\varrho((r'_i)_{i\in I})\otimes x+r'\varrho((r_i)_{i\in I})\otimes x+xr'\varrho((r'_i)_{i\in I})\otimes x)  \\
   & = &(r\varrho((r_i)_{i\in I}),\varrho((r_i)_{i\in I})\otimes r'x+\varrho((r'_i)_{i\in I})\otimes(r+xr')x)  \\
   & = & (r,r+xr')\ast(\varrho((r_i)_{i\in I},\varrho((r'_i)_{i\in I}\otimes x) \\
   & = & (r,r+xr')\ast \widetilde{\varrho}((r_i,r_i+r'_ix))_{i\in I})
\end{eqnarray*}
Consequently, $\widetilde{\rho}$ is an $R\bowtie xR$-morphism from
$\widehat{M}$ into $(R\bowtie xR)^{(I)}$ and $\widehat{\varrho}$ is
an $R\bowtie xR$-morphism from $(R\bowtie xR)^{(I)}$ into
$\widehat{M}$.\\ Moreover, if $(\overline{m},\overline{m'}\otimes
x)\in \overline{M}\oplus(\overline{M}\otimes_RxR)$, we have
\begin{eqnarray*}
  \widehat{\varrho}\circ\widehat{\rho}(\overline{m},\overline{m'}\otimes x) &=& \widehat{\varrho}((\rho_i(\overline{m}),\rho_i(\overline{m})+x\rho_i(\overline{m'}))_{i\in
I}) \\
   &=& (\varrho((\rho_i(\overline{m})_{i\in I})),\varrho((\rho_i(\overline{m'})_{i\in I})\otimes x) \\
   &=&  (\varrho(\rho(\overline{m})),\varrho(\rho(\overline{m'}))\otimes x)\\
   &=& (\overline{m}, \overline{m'}\otimes x)
\end{eqnarray*}
 Thus, $\widehat{M}$ is a projective $R\bowtie xR$-module. Hence, $M$ is also
 projective as an $R\bowtie xR$-module (recall that
 $\widehat{M}\cong_{R\bowtie xR}M$ by the $R\bowtie xR$-isomorphism $\beta$). Hence, the case $n=0$ holds.\\
Now,  suppose that $0<\pd_R(M\otimes_{R\bowtie xR}R)= n$. Clearly
$\pd_{R\bowtie xR}(M)>0$. Thus,  consider a short exact sequence of
$R\bowtie xR$-modules:
$$(\mp)\quad 0\longrightarrow K \longrightarrow P \longrightarrow M
\longrightarrow 0$$ where $P$ is projective and $\pd_{R\bowtie
xR}(K)= \pd_{R\bowtie xR}(M)-1$. Since $\Tor_{R\bowtie xR}^i(M,R)$
for every $i>0$, if we tensor $(\mp)$ with $-\otimes_{R\bowtie
xR}R$, we obtain the short exact sequence of $R$-modules:
$$0\longrightarrow K\otimes_{R\bowtie xR}R\longrightarrow
P\otimes_{R\bowtie xR}R\longrightarrow M\otimes_{R\bowtie
xR}R\longrightarrow 0$$ and moreover, we have $\Tor_{R\bowtie
xR}^i(K,R)$ for every $i>0$. Thus, by hypothesis, $$\pd_{R\bowtie
xR}(M)=\pd_{R\bowtie xR}(K)+1\leq \pd_R(K\otimes_{R\bowtie
xR}R)+1=\pd_R(M\otimes_{R\bowtie xR}R)$$ as desired.
\end{proof}

Recall that the finitistic projective dimension of a ring $R$ is

$$\FPD(R)=\sup\{\pd_R(M)\mid \text{M  is  an  R-module and
$\pd_R(M)<\infty$}\}.$$

\begin{lemma}\label{FPD}
Let $R$ be a ring which contains a regular element $x$. Then,
$\FPD(R\bowtie xR)=\FPD(R)$.
\end{lemma}
\begin{proof}
First, we claim that $\FPD(R\bowtie xR)\leq \FPD(R)$. Obviously,
if $\FPD(R)=\infty$ the inequality holds. Hence we may assume that
$\FPD(R)<\infty$. Let $M$ be an $R\bowtie xR$-module with finite
projective dimension. From the shorts exact $R\bowtie
xR$-sequences:
$$(\star)\quad 0\longrightarrow 0\times
xR\stackrel{\iota}\longrightarrow R\bowtie
xR\stackrel{\varepsilon}\longrightarrow R \longrightarrow 0$$
where $i$ is the injection and $\varepsilon(r,r+r'x)=r$, and
$$(\star\star) \quad 0 \longrightarrow R\stackrel{\mu}\longrightarrow
R\bowtie xR \stackrel{\nu}\longrightarrow 0\times xR\longrightarrow
0$$ where $\mu(r)=(rx,0)$ and $\nu(r,r+xr')=(0,(r+r'x)x),$\\
we deduce that, for every $i>0$,  $$\Tor^i_{R\bowtie
xR}(R,M)=\Tor^{i+1}_{R\bowtie xR}(0\times
xR,M)=...=\Tor^{i+m}_{R\bowtie xR}(R,M)=0$$ for some integer $m>n$.
Thus, Using Lemma \ref{passage}, $\pd_R(M\otimes_{R\bowtie xR}R)=
\pd_{R\bowtie xR}(M)<\infty$. So, we  conclude that $\pd_{R\bowtie
xR}(M)=\pd_R(M\otimes_{R\bowtie xR}R)\leq \FPD(R)$.
Consequently, $\FPD(R\bowtie xR)\leq \FPD(R)$, as desired.\\
Secondly, we claim that $\FPD(R)\leq \FPD(R\bowtie xR) $. So, let
$M$ be an $R$-module with finite projective dimension. Since
$R\bowtie xR\cong_R R\oplus xR\cong_RR^2$ (since $x$ is regular
element) it is easy to see that $\pd_{R\bowtie
xR}(M\otimes_RR\bowtie xR)\leq \pd_R(M)<\infty$. Thus,
$\pd_{R\bowtie xR}(M\otimes_RR\bowtie xR)\leq \FPD(R\bowtie xR)$.
And we have:
$$\pd_R(M)=\pd_R(M\otimes_RR\bowtie xR)\leq \pd_{R\bowtie xR}(M\otimes_RR\bowtie xR)\leq
\FPD(R\bowtie xR)$$

 Consequently, $\FPD(R)\leq \FPD(R\bowtie xR)$, as desired.\\ From the first and second claim,
 we conclude the desired equality.
\end{proof}

In \cite{Ding}, Ding and Mao introduced the cotorsion dimension of
modules and rings, which are defined as follows:
\begin{definition}
Let $R$ be a ring.
\begin{enumerate}
  \item The cotorsion dimension of an $R$-module $M$, denoted by $\cd_R(M)$,
is the least positive integer $n$ for which $\Ext^{n+1}_R(F,M)= 0$
for all flat $R$-modules $F$.
  \item The global cotorsion dimension of $R$,
denoted by $\cotd(R)$, is defined as the supremum of the cotorsion
dimensions of R-modules.
\end{enumerate}

\end{definition}
 In \cite[Theorem 7.2.5(1)]{Ding}, the authors prove that we have, for a ring $R$ and a positive integer $n$,
 $\cot.D(R) \leq n$ if and only if, every flat $R$-module $F$ has a projective dimension less or equal than
$n$. In \cite[Theorem 2.1]{bennis Mahdou 4}, we find for any
coherent ring the relation $\cotd(R) \leq \Ggldim(R) \leq \wGgldim(R)
+ \cotd(R)$.

\begin{lemma}\label{cot}
Let $R$ be a ring which contains a regular element  $x$. Then,
$\cotd(R\bowtie xR)=\cotd(R)$.
\end{lemma}

\begin{proof}
First, we claim that $\cotd(R\bowtie xR)\leq \cotd(R)$. From,
\cite[Theorem 7.2.5(1)]{Ding}, it remains to prove that the
projective dimension of every flat $R\bowtie xR$ is less or equal
than  $\cotd(R)$. So, let $F$ be a flat $R\bowtie xR$-module. Then,
from Lemma \ref{passage}, $\pd_{R\bowtie
xR}(F)=\pd_R(F\otimes_{R\bowtie xR}R)$. On the other hand,
$F\otimes_{R\bowtie xR}R$ is a flat $R$-module. Thus, by
\cite[Theorem 7.2.5(1)]{Ding}, $\pd_R(F\otimes_{R\bowtie xR}R)\leq
\cotd(R)$. Therefore, $\pd_{R\bowtie xR}(F)\leq \cotd(R)$.
Consequently,
$\cotd(R\bowtie xR)\leq \cotd(R)$, as desired.\\
Secondly, we claim that $\cotd(R)\leq \cotd(R\bowtie xR)$. Similarly to
the first claim, we have to prove that every flat $R$-module has a
projective dimension less or equal than $\cotd(R\bowtie xR)$. So,
let $F$ be a flat $R$-module. By \cite[Theorem 7.2.5(1)]{Ding},
$\pd_{R\bowtie xR}(F\otimes_RR\bowtie xR)\leq \cotd(R\bowtie xR)$.
Moreover, since $R\bowtie xR\cong_RR^2$ (recall that $x$ is a
regular element of $R$), we have:
$$\pd_R(F)=\pd_R(F\otimes_R R\bowtie xR)\leq \pd_{R\bowtie
xR}(F\otimes_RR\bowtie xR)\leq \cotd(R\bowtie xR)$$ Consequently,
$\cotd(R)\leq \cotd(R\bowtie xR)$, as desired. Hence, by the first
and second claim we have the desired result.
\end{proof}

\begin{proof}[Proof of Theorem \ref{main result2}]
If $R$ is coherent, $R\bowtie xR$ is also coherent (\cite[Theorem
3.1]{chhitimahdou}).\\

 $(1)$ Follows immediately from Lemmas  \ref{n-FC} and  \ref{FP-id}.\\

$(2)$ The proof of this result is the same as the one of
\cite[Proposition 2.5]{bennis Mahdou 4}. For completeness, we give a
proof here.\\
First assume that $\Ggldim(R)$ is finite. By Lemma \ref{FPD} and
\cite[Theorem 2.28]{Holm}, $\FPD(R\bowtie xR)=\FPD(R)=\Ggldim(R)$ is
finite. On the other hand, from \cite[Theorem 2.1]{bennis Mahdou 4}
and Lemma \ref{cot}, we have $\cotd(R\bowtie xR)=\cotd(R)\leq \Ggldim(R)$ is
finite. And, by  \cite[Corollary 1.2]{Bennis and Mahdou2} and $(1)$
above, $\wGgldim(R\bowtie xR)=\wGgldim(R)\leq \Ggldim(R)$ is finite.
Consequently, from \cite[Theorem 2.1]{bennis Mahdou 4},
$\Ggldim(R\bowtie xR)\leq \wGgldim(R\bowtie xR)+\cotd(R\bowtie xR)$ is
finite. Therefore, from \cite[Theorem 2.28]{Holm}, $\Ggldim(R\bowtie
xR)=\FPD(R\bowtie xR)=\Ggldim(R)$.\\
Similarly, we show that $\Ggldim(R\bowtie xR)=\Ggldim(R)$ when
$\Ggldim(R\bowtie xR)$ is finite, and this gives the desired result.
\end{proof}

From Theorem \ref{main result2} and \cite[Corollary 1.2]{Bennis
and Mahdou2} we see that for every coherent ring $R$ with finite
classical homological dimensions,  the \emph{amalgamated
duplication} of $R$ along a principal regular  ideal $I=xR$ transfers
$R$ to a ring $R\bowtie xR$ with infinite classical homological
dimensions but with
 finite Gorenstein dimensions. Namely, $\Ggldim(R\bowtie
 xR)=\gldim(R)$ and $\wGgldim(R)=\wdim(R)$.\\

 Now, we are able to construct a new
class of reduced Noetherian rings with finite Gorenstein global
dimensions and infinite weak dimensions.

\begin{example}
 Let $R$ be a Noetherian domain of global dimension equal to
$n$, $x$ be a non invertible regular element of $R$ and let $S:=R
\bowtie xR$. Then:
\begin{enumerate}
  \item $\wGgldim(S)\; (=\Ggldim(S)) =n$ (by Theorem \ref{main result2} and \cite[Theorem 12.3.1]{Enocks and janda}).
  \item $\wdim(S)\; (=\gldim(S)) =\infty$ and $S$ is a Noetherian reduced ring (by
  \cite[Proposition 2.1]{D'Anna}).
\end{enumerate}
\end{example}
In the following example, we construct a new family of
non-Noetherian coherent rings $\{S_n\}_{n\geq 1}$ such that
$\Ggldim(S_n)=n+1$ and $\wdim(S_n)=\infty$ for each $n>0$.

\begin{example} Let $n > 0$ be an integers and let $R_n=R[X_1,X_2,...,X_n]$ be the polynomial ring
in $n$ indeterminates over a non-Noetherian hereditary ring $R$.
Let $S_n:=R_n\bowtie X_1R_n$. Then, $S_n$ is a non Noetherian
coherent ring with $\Ggldim(S_n)=n+1$ and $\wdim(S_n)=\infty$.
\end{example}

\begin{proof}
From \cite[Theorem 7.3.1]{Glaz}, $R_n$ is coherent for every
$n>0$. And by Hilbert's syzygy Theorem, $\gldim(R_n)=n+1$.
Therefore, Theorem \ref{main result2} and \cite[Theorem
2.1]{Bennis and Mahdou3} imply that $\Ggldim(S_n)=n+1$.
\end{proof}

Also, Theorem \ref{main result2} allows us to construct a new
family of non-Noetherian coherent rings $\{R_n\}_{n\geq 1}$ such
that $\wGgldim(R_n)=n+1$, $\Ggldim(R_n) > n+1$ and
$\wdim(R_n)=\infty$ for all $n\geq 0$, as follows:

\begin{example}
Let $R$ be a non-semisimple quasi-Frobenius ring and let $S$ be a
non-Noetherian Pr\^ufer domain. Let $S_n:=R\times S[X_1,..,X_n]$
for every $n>0$, $S_0=R\times S$, $a$ a non-unit element of $S$,
$x=(1,a)$ and $R_n:=S_n\bowtie xS_n$ for every $n\geq 0$. Then,
$\wGgldim(R_n)=n+1$, $\Ggldim(R_n)>n+1$ and $\wdim(R_n)=\infty$.
\end{example}

\begin{proof}
From \cite[Example 3.8]{Bennis and Mahdou3}, $S_n$ is a
non-Noetherian coherent ring, $\wGgldim(S_n)=n+1$ and
$\Ggldim(S_n)>n+1$ for every $n$. See also that $x$ is a non-unit
regular element of $S_n$ for every $n\geq 0$. Thus, by
\cite[Theorem 3.1]{chhitimahdou}, \cite[Proposition 2.1]{D'Anna}
and Theorem \ref{main result2}, $S_n$ is non-Noetherian coherent,
$\wGgldim(R_n)=n+1$ and $\Ggldim(R_n)>n+1$ and
$\wdim(R_n)=\infty$, for every $n$, as desired.
\end{proof}

\bibliographystyle{amsplain}

\end{document}